\documentclass[leqno,12pt]{article}
\usepackage{amssymb,amsfonts}
\usepackage{amsmath,latexsym}

\newcommand{\bepr}{{\em Proof} } 
\newcommand{\enpr}{\hfill \rule{.5em}{.5em}}

\newcommand{\R}{{\mathbb R}}

\newtheorem{prop}{Proposition}[section] 
\newtheorem{thm}{Theorem}[section] 
\newtheorem{lemma}{Lemma}[section] 

\newtheorem{cor}{Corollary}[section]

\begin{document}

\title{The role of the Hilbert metric in a class of \\ singular elliptic boundary value problem in convex domains}

\author{Denis Serre \\ UMPA, UMR CNRS--ENS Lyon \# 5669 \\ \'Ecole Normale Sup\'erieure de Lyon\thanks{46, all\'ee d'Italie, F--69364 Lyon, cedex 07}}

\date{\today}

\maketitle

\begin{abstract} In a recent paper \cite{Ser_log}, we were led to consider a distance over a bounded open convex domain. It turns out to be the so-called {\em Thompson metric}, which is equivalent to the {\em Hilbert} metric. It plays a key role in the analysis of existence and uniqueness of solutions to a class of elliptic boundary-value problems that are singular at the boundary.
\end{abstract}

\section{Introduction}

Let $\Omega$ be a connected open set in $\R^n$. If $p\in\R^n$, we denote $|p|$ its usual Euclidian norm. The class of boundary value problems that we are interested in is
\begin{eqnarray}
\label{eq:PDE}
{\rm div}(a(|\nabla w|)\nabla w)+\frac{F(|\nabla w|)}w & = 0 & \hbox{ in }\Omega, \\
\label{eq:posit}
w & > 0 & \hbox{ in }\Omega, \\
\label{eq:DBC}
w & = 0 & \hbox{ on }\partial\Omega.
\end{eqnarray}
Hereabove, $a$ is a smooth numerical even function, which satisfies the requirements for ellipticity:
\begin{equation}
\label{eq:ellipass}
a(r)>0,\qquad a(r)+ra'(r)>0,\qquad \forall r\ge0.
\end{equation}
We warn the reader that we do not assume {\em a priori} a uniform ellipticity~; it may happen that the ratio
$$\frac{a(r)}{a(r)+ra'(r)}$$
tends either to $0$ or to $+\infty$ as $r\rightarrow+\infty$. For instance, we allow the principal part to be the minimal surface operator, where $a(r)=(1+r^2)^{-1/2}$, for which 
$$\frac{a(r)}{a(r)+ra'(r)}=1+r^2\rightarrow+\infty.$$

We suppose that $F$ is a smooth, non-negative function, and that 
\begin{equation}\label{eq:Fass}
F(0)>0.
\end{equation}
The lower order term in (\ref{eq:PDE}) therefore becomes singular at the boundary, where $w$ vanishes. 

\paragraph{Notations.} In the sequel, we denote $b(r)=ra(r)$, so that $a+ra'=b'$. We define a strictly increasing function
$$G(r)=\int_0^r\frac{sb'(s)}{F(s)}\,ds.$$
The inverse $G^{-1}$ will be denoted $H$.

\subsection{Data}

At first glance, it may look strange that neither the equation, nor the boundary contain some explicit data~; both are ``homogeneous''. Our data is nothing but the domain itself. 
The assumption about $\Omega$ meets that in other works on non-uniformly elliptic BVPs: it is a {\em bounded convex} domain in $\R^n$.

\subsection{Motivations}

We came to this class of problems through the analysis of the two-dimensional Riemann problem for the Euler system of a compressible flow, when the gas obeys the so-called Chaplygin equation of state. This problem can be recast as 
\begin{equation}\label{eq:RP}
{\rm div}\frac{\nabla w}{\sqrt{1+|\nabla w|^2}}+\frac2{w\sqrt{1+|\nabla w|^2}}=0,
\end{equation}
which is (\ref{eq:PDE}) with
$$n=2,\qquad a(r)=\frac1{\sqrt{1+r^2}}\,,\qquad F=2a.$$
We first proved the existence and uniqueness (see \cite{Ser_Chap}) for this problem whenever $\Omega$ is uniformly convex, in the sense that the curvature is bounded away from zero along the boundary. Later on, we removed the uniform assumption and proved the existence for every convex bounded planar domain \cite{Ser_log}. This improvement involves an interior Lipschitz estimate of $\log w$ in terms of a special metric over $\Omega$, for which the boundary is a horizon. We shall show below that this distance is nothing but the Hilbert metric $d_H$, giving meanwhile a new and rather simple proof of the triangle inequality.

It turns out that the very same BVP also governs those graphs $x_3=w(x_1,x_2)$ that are complete minimal surfaces in the $3$-dimensional hyperbolic space ${\mathbb H}_3$, the upper half-space in $\R^3$, equipped with the metric 
$$ds^2=\frac{dx_1^2+dx_2^2+dx_3^2}{x_3^2}\,,$$
of constant negative curvature. The existence of such minimal surfaces was studied by Anderson \cite{And} in the parametric and the non-parametric cases, the latter involving the graph of $w$. The non-parametric part of Anderson's paper is however incomplete, in that the author  contents himself to establish $L^\infty$-bounds (by below and above) and claims that it automatically implies regularity estimates in the interior. This claim is not true because the principal part of the PDE, the operator for minimal surfaces, is not {\em a priori} uniformly elliptic. Uniform ellipticity requires the knowledge of a prior Lipschitz estimate, which can not be overlooked. The same flaw occurs in Lin's paper \cite{Lin}.

We point out that in both of these motivations, the convexity of $\Omega$ is a necessary condition for existence (and therefore a necessary and sufficient one). In the Chaplygin Riemann problem, it is guaranted by the analysis of the propagation of shock waves. If $\Omega$ is not convex, a complete minimal surface in ${\mathbb H}_3$, asymptotic to $\partial\Omega$, exists\footnote{This is the parametric part, by far the main one,  of Anderson's paper, on which we have no doubt at all.} as a current \cite{And}, but it is not a graph over $\Omega$.

\subsection{Content of this paper}

We start by showing in Section \ref{s:dist} the equality between our (not so) new distance and the Hilbert metric in $\Omega$.

Then we turn towards  the class of BVPs (\ref{eq:PDE},\ref{eq:posit},\ref{eq:DBC}).
We show  that essentially the same strategy as the one designed in \cite{Ser_log} works out under the rather mild assumption that
\begin{equation}
\label{eq:assump}
\int^{+\infty}e^{-G(s)}\,\frac{b'(s)}{F(s)}\,ds<\infty.
\end{equation}
Our main result is therefore
\begin{thm}\label{th:main}
Let $a,F$ be even smooth functions, satisfying (\ref{eq:ellipass},\ref{eq:Fass},\ref{eq:assump}). Then, for every bounded convex domain $\Omega\subset\R^d$, there exists one and only one function
$$w\in{\cal C}(\overline{\Omega})\cap{\cal C}^\infty(\Omega)$$
solving the BVP (\ref{eq:PDE},\ref{eq:posit},\ref{eq:DBC}).

In addition, $\log w$ is Lipschitz, with constant $1$, with respect to the Hilbert metric $d_H$.
\end{thm}
Of course, if $a$ or $F$ has only finite regularity, then $w$ has only finite regularity.

\subsection*{Acknowledgement}

I am indebted to Ludovic Marquis for useful comments about the Hilbert metric and its variants.

\section{A distance over a bounded convex domain}\label{s:dist}

Let $\Omega$ is a non-void,  bounded convex open domain in $\R^n$. 
Given two points $p,q\in\Omega$, $\Omega-p$ contains a ball centered at the origin and is therefore absorbing. Thus there exists some $\lambda>0$ such that $\Omega-q\subset\lambda(\Omega-p)$. If $\mu>\lambda$, then also $\Omega-q\subset\mu(\Omega-p)$, by convexity. Likewise, the infimum $m(p,q)$ of all such numbers satisfies the same inclusion, by continuity. Hence the set of these numbers is of the form $[m(p,q),+\infty)$. Considering the volumes, we have
$$|\Omega|=|\Omega-q|\le m(p,q)^n|\Omega-p|=m(p,q)^n|\Omega|,$$
which implies 
\begin{equation}
\label{eq:mgeun}
m(p,q)\ge1.
\end{equation}
The equality in (\ref{eq:mgeun}) stands only if 
$$\Omega-q=\Omega-p,$$
that is if $p=q$.

If $r\in\Omega$ is a third point, then
$$\Omega-r\subset m(q,r)(\Omega-q)\subset m(q,r)m(p,q)(\Omega-p)$$
and therefore
$$m(p,r)\le m(q,r)m(p,q).$$
All this shows that the logarithm of $m$ is a non-negative function over $\Omega\times\Omega$, which vanishes only along the diagonal and satisfies the triangle inequality. In other words, the function
$$d_\Omega(p,q)=\log m(p,q)+\log m(q,p)$$
is a {\em distance} over $\Omega$. In our paper \cite{Ser_log}, we used the equivalent metric 
$$d_\Omega'(p,q)=\max\{\log m(p,q),\log m(q,p)\}.$$

\bigskip

We prove here that $d_\Omega$ is nothing but the Hilbert distance $d_H$ over $\Omega$ (see \cite{Hil}). Let us recall the definition of the latter. If $p,q\in\Omega$, let $r,s$ be the intersection points of the line $L$ passing through $p$ and $q$, with the boundary $\partial\Omega$~; we label the points so that $r,p,q,s$ are in this order along $L$. Then $d_H$ is the logarithm of a cross-ratio~:
$$d_H(p,q)=\log\frac{\overline{rq}\cdot\overline{ps}}{\overline{rp}\cdot\overline{qs}}\,.$$
Our first result therefore reads
\begin{prop}\label{p:idem}
For every  non-void,  bounded convex open domain $\Omega\subset\R^n$, one has
\begin{equation}\label{eq:idem}
d_\Omega\equiv d_H.
\end{equation}
\end{prop}

\bigskip

The equality (\ref{eq:idem}) follows immediately from the
\begin{lemma}\label{l:mpq}
With the notation above, there holds
\begin{equation}\label{eq:mla}
m(p,q)=\frac{\overline{rq}}{\overline{rp}}\,.
\end{equation}
\end{lemma}

\bepr

Let us define $\theta=\frac1u\in(0,1)$. Because $r$ is a boundary point of $\Omega$, and $\Omega$ is convex open, we have $(1-\theta)\Omega+\theta r\subset\Omega$. This is equivalent to writing 
$$\Omega-q\subset\frac{u}{u-1}(\Omega-p),$$
which is the inequality $\le$ in (\ref{eq:mla}).

Conversely, suppose that $\Omega-q\subset\lambda(\Omega-p)$. Then we have $\overline{\Omega}-q\subset\lambda(\overline{\Omega}-p)$ and therefore $r-q\in\lambda(\overline{\Omega}-p)$. This amounts to writing
$$\left(1+\frac u\lambda\right)p-\frac u\lambda\in\overline{\Omega},$$
but this is equivalent to
$$v\le1+\frac u\lambda\le u.$$
The second inequality gives $\lambda\ge\frac u{u-1}$\,. This implies the inequality $\ge$ in (\ref{eq:mla}).

\enpr

\bigskip

\paragraph{Remarks.} This characterization of the Hilbert metric is related to the construction of the Hilbert projective metric over the cone 
$$C=\{(t,tx)\,|\, t>0\hbox{ and } x\in\Omega\},$$
see \cite{Wal}. 
It provides a much simpler proof of the triangle inequality than the original one. For the classical proof, which involves projetive geometry, see the introductory article in {\em Image des Math\'ematiques} \cite{IdM}.
Lemma \ref{l:mpq} also implies that $d_\Omega'$ is identical to the Thompson metric in $\Omega$.

\section{The strategy for existence and uniqueness to the BVP}

Our first observation is that the PDE (\ref{eq:PDE}) is of the quasilinear form 
$$\sum_{i,j}a_{ij}(\nabla w)\partial_i\partial_jw+N(w,\nabla w)=0,$$
where the principal part is elliptic,
$$c(p)|\xi|^2\le\sum_{i,j}a_{ij}(p)\xi_i\xi_j\le C(p)|\xi|^2,\qquad 0<c(p)<C(p)<\infty.$$
The coefficients $a_{ij}$  involve the gradient of $w$, but not $w$ itself. Finally the lower order term $N$ is non-increasing in $w$. Therefore the PDE satisfies the maximum principle (MP).

The MP allows us to compare a sub-solution and a super-solution. A locally Lipschitz function $u:\Omega\rightarrow(0,+\infty)$ is a sub-solution of (\ref{eq:PDE}) if it satisfies, in the distributional sense,
\begin{equation}
\label{eq:subs}
{\rm div}(a(|\nabla u|)\nabla u)+\frac{F(|\nabla u|)}u\ge0.
\end{equation}
It is a super-solution if it satisfies the opposite inequality $\le$ in (\ref{eq:subs}). If in addition $u$ is continuous over $\overline\Omega$, we say that $u$ is a super-solution of the BVP if it is a super-solution of (\ref{eq:PDE}), and it satisfies $u\ge0$ on $\partial\Omega$. It is a sub-solution if it satisfies (\ref{eq:subs}), and $u\le0$ over $\partial\Omega$ (but then this means $u\equiv0$ on the boundary, because $u>0$ in the interior).

If $u$ and $v$ are a sub-solution and a super-solution respectively, of the BVP in some domain $\Omega$, then $u\le v$ in $\Omega$. In particular, if $w$ is a solution in $\Omega$, then $u\le w\le v$. This immediately implies the uniqueness part of Theorem \ref{th:main}.

The method for existence is based on the one hand on a continuation argument, described in Section \ref{s:exist}, and on the other hand on {\em a priori} estimates. The latter must be robust enough to allow us to pass to the limit in a sequence of solutions. To ensure the boundary condition, we shall use sub- and super-solution respectively to construct barrier functions $w_\pm>0$, with $w_\pm\equiv0$ on the boundary. The fact that $w$ is clamped between $w_-$ and $w_+$ implies the boundary condition. It also ensures that $w$ is positive and bounded in $\Omega$. In order to pass to the limit in the PDE, we need a precompactness property of $\nabla w$ in $L^\infty_{loc}(\Omega)$. This will be given by a ${\cal C}^2_{loc}$-regularity estimate and the Ascoli--Arzela theorem. The regularity  is a well-known fact (see Gilbarg \& Trudinger \cite{GT}) whenever the operator 
$$L(p)=\sum_{i,j}a_{ij}(p)\partial_i\partial_j$$
is uniformly elliptic. Since we have not assumed the latter property, it must come as a consequence of the fact that $p=\nabla w$ takes its values in a compact subset of $\R^n$. In other words, we need an {\em a priori} estimate  of $\nabla w$ in $L^\infty_{loc}(\Omega)$.

We summarize below the tasks we are going to address:
\begin{itemize}
\item Construct a finite upper bound $w_+$ of $w$, continuous up to the boundary, where it satisfies $w_+\equiv0$.
\item Construct a lower bound $w_->0$ of $w$, continuous up to the boundary, where it satisfies $w_-\equiv0$.
\item Find a Lipschitz estimate of $w$ in $\Omega$. This estimate may deteriorate near the boundary, but it must be uniform on every compact subset of $\Omega$. This is where the Hilbert metric is at stake.
\item Make all these estimates uniform with respect to some approximation.
\end{itemize}
Of course, the only tool at our disposal is the maximum principle.

\section{The barrier functions}

We shall use repeatedly the fact that the PDE (\ref{eq:PDE}) is invariant under a scaling: if $z$ is a solution in some domain $\omega$, then the function $z^\mu(x):=\frac1\mu z(\mu x)$ is again a solution, in $\frac1\mu\,\omega$.

\subsection{The upper barrier}\label{ss:upper}

We write our convex domain as the intersection of slabs
$$\Omega=\bigcap_{\nu\in S^{n-1}}\Pi_\nu,\qquad\pi_\nu=\{x\in\R^n\,|\,\alpha_-(\nu)< x\cdot\nu<\alpha_+(\nu)\},$$
where we have of course $\alpha_\pm(-\nu)=-\alpha_\mp(\nu)$. Notice that $\nu\mapsto\alpha_\pm$ is continuous.

Our upper bound will be given as the infimum of super-solutions. The building block is the solution of the BVP in the interval $(0,1)$~:
\begin{lemma}[$1$-D case.]\label{l:W}
When $n=1$ and the domain is $(0,1)$, then the BVP admits a unique solution $W$.
\end{lemma}
We infer that the BVP in a slab $\Pi_\nu$ admits a solution, namely
$$W_\nu(x)=(\alpha_+(\nu)-\alpha_-(\nu))W\left(\frac{x\cdot\nu-\alpha_-(\nu)}{\alpha_+(\nu)-\alpha_-(\nu)}\right).$$
Because $\Omega\subset\Pi_\nu$ and $W_\nu$ is non-negative, in particular along $\partial\Omega$, its restriction to $\Omega$ is a super-solution of the BVP in $\Omega$. Therefore the solution $w$ satisfies $w\le W_\nu$. This yields to our upper-bound,
$$w(x)\le w_+(x)=\inf_{\nu\in S^{n-1}}W_\nu(x).$$
The continuity of $\alpha_\pm$, plus the uniform continuity of $W$, imply that $w_+$ is continuous over $\overline\Omega$. We point out that, because every $y\in\partial\Omega$ is a boundary point of some $\Pi_\nu$, $w_+$ vanishes on the boundary.

\bigskip

\bepr 

We already know the uniqueness. Using the reflexion $x\leftrightarrow 1-x$, we infer that $W$ must be even: $W(1-x)=W(x)$. We anticipate that $W$ is monotonous over $(0,\frac12)$ and write the PDE, now an ODE as
$$(b(W'))'+\frac{F(W')}W=0,\qquad W(0)=W'(\frac12)=0,$$
where $b(r):=ra(r)$. We recall that $b'>0$, from ellipticity.

Let us define $z=W'\circ W^{-1}$. Using $W'=z(W)$, we transform the ODE into 
$$zb'(z)\frac{dz}{dW}+\frac{F(z)}W=0.$$
The latter ODE amounts to writing $G(z)+\log W={\rm cst}$, from which we obtain $z=H(\log\frac cW)$ for some integrating factor $c\in\R$.

Let us make temporarily the choice that $c=1$ and consider a maximal solution of  the autonomous ODE $W'=H(-\log W)$. We have
$$\frac{dW}{H(-\log W)}=dx.$$
Because of (\ref{eq:assump}), 
we have
$$\int_0\frac{ds}{H(-\log s)}\,<\infty.$$
Therefore there exists a unique solution $W_0$ of the Cauchy problem $$W_0'=H(-\log W_0),\qquad W_0(0)=0.$$ 
This $W_0$ is increasing. Since the integral
$$\int_0e^{-G(s)}\frac{b'(s)}{F(s)}\,ds$$
is converging, we have
$$\int^{e^{-G(0)}}\frac{ds}{H(-\log s)}\,<\infty$$
and therefore $W_0$ reaches the value $e^{-G(0)}$ at some finite $\bar x>0$. Then $W_0'(\bar x)=0$. Extending $W_0$ it by parity, we obtain a solution of the BVP in the interval $(0,2\bar x)$. Then 
$$W(t)=\frac 1{2\bar x}W_0(2\bar xt)$$
defines the solution of the BVP over $(0,1)$.

\enpr

\subsection{The lower barrier}

The construction of the lower barrier does not make use of the convexity. We begin with a building block:
\begin{lemma}
There exists an $\epsilon>0$ such that the function $Z(x)=\frac\epsilon2(1-|x|^2)$ be a sub-solution of the BVP in the unit ball $B(0;1)$.
\end{lemma}

\bepr

Since $Z$ is positive in the ball, it suffices to check that  $Z$ satifies (\ref{eq:subs}). This inequality writes
$$\frac12(\epsilon^2-t^2)(b(t)+(d-1)a(t))\le F(t),\qquad\forall t\in[0,\epsilon].$$
Because $a$ and $b$ are non-negative, it is enough to have
$$\frac{\epsilon^2}2(b(t)+(d-1)a(t))\le F(t),\qquad\forall t\in[0,\epsilon].$$
Let $A$ and $B$ be the upper bounds of $a$ and $b$ over $[0,1]$ respectively. If $\epsilon<1$, it is enough to have
$$\frac{\epsilon^2}2(B+(d-1)A)\le \sup_{t\in[0,\epsilon]}F(t),$$
which is obviously true for $\epsilon>0$ small enough.

\enpr

\bigskip

By translation and scaling, we inherit a sub-solution of the BVP in any ball $B(x_0;\rho)$~:
$$Z_{x_0,\rho}(x)=\rho Z\left(\frac{x-x_0}\rho\right).$$
If $B(x_0;\rho)$ is contained in $\Omega$, then $w$ is a super-solution for the BVP in this ball, and we infer $w\ge Z_{x_0,\rho}$. This leads us to our lower barrier function 
$$w_-(x)=\sup\{Z_{x_0,\rho}\,|\,B(x_0;\rho)\subset\Omega\}.$$
We point out that $w_-$ is continuous over $\overline\Omega$ and is positive in the interior.

\section{The Lipschitz estimate}

The main ingredient is the
\begin{lemma}\label{l:LipH}
The solution of the BVP (\ref{eq:PDE},\ref{eq:posit},\ref{eq:DBC}) in a bounded convex open domain $\Omega$ satisfies, if it exists
\begin{equation}\label{eq:unL}
|\log w(q)-\log w(p)|\le \max\{\log m(p,q),\log m(q,p)\},\qquad\forall p,q\in\Omega.
\end{equation}
Consequently, $\log w$ is Lipschitz with constant at most $1$, with respect to the Hilbert metric.
\end{lemma}

Because the restriction to the Hilbert metric to a compact subset $K\subset\Omega$ is equivalent to the Euclidian distance, we infer a Lipschitz estimate in the classical sense, away from the boundary. Because $\min_K w_->0$ and $w_+$ is bounded, this transfers into a local Lipchitz estimate of $w$:
\begin{cor}
For every compact subset $K\subset\Omega$, the restriction $w|_K$ enjoys an {\em a priori} estimate in the Lipschitz semi-norm $\sup_K|\nabla w|$.
\end{cor}

\bigskip

\bepr

Given $p,q\in\Omega$, the function
$$x\mapsto m(p,q)w\left(\frac{x}{m(p,q)}+p\right)$$
is the solution of the BVP in the domain $m(p,q)(\Omega-p)$. Since the latter contains $\Omega-q$, it is also a super-solution in the domain $\Omega-q$. It is therefore larger than or equal to the solution $w(x+q)$ in the latter:
$$w(x+q)\le m(p,q)w\left(\frac{x}{m(p,q)}+p\right),\qquad\forall x\in\Omega-q.$$
Setting $x=0$ in the inequality above, we derive
$$w(q)\le m(p,q)w(p).$$
Exchanging the roles of $p$ and $q$, we also have $w(x)\le m(q,p)w(y)$, whence (\ref{eq:unL}).

\enpr

\subsection{The best Lipschitz constant}

Lemma \ref{l:LipH} provides an upper bound for the Lipschitz constant of $\log w$ with respect to the Hilbert metric:
$$c_\Omega:=\sup_{x\ne y}\frac{|\log w(y)-\log w(x)|}{d_H(x,y)}\,\le1.$$
We may wander whether this bound is accurate or not. Remark that if $O$ is a boundary point and $L$ is a ray emanating from $O$ in $\Omega$, then the restriction of $d_H$ to $L$ is logarithmic, in the sense that if $x,y\in L$, then 
$$d_H(x,y)=|\log t_y+\log(T_L-t_x)-\log t_x-\log(T_L-t_y)|\sim|\log t_y-\log t_x|,$$
where $t$ is the affine coordinate along $L$ with origin $O$, and $T_L$ is the coordinate or the other intersection point of $L$ with $\partial\Omega$. If the solution admits a H\"older singularity at a boundary point, of exponent $\alpha\in(0,1]$, we deduce that $c_\Omega\ge\alpha$. 

One remarquable application of this principle is the following
\begin{prop}\label{p:coneL}
Let the origin be a conical point of $\partial\Omega$, and denote ${\cal C}$ the tangent cone at $0$. Suppose that  the BVP is solvable in the cone ${\cal C}$, with a solution $V(x)=|x|v\left(\frac x{|x|}\right)$. Then the solution in $\Omega$ is asymptotic to $V$ as $x\rightarrow0$. In particular,
$$c_\Omega=1.$$
\end{prop}
The fact that $V$ is homogeneous of degree one is a consequence of the scaling invariance of ${\cal C}$ (the conic property) and the expected uniqueness.

\bigskip

\bepr

Let $w$ be the solution of the BVP in $\Omega$, and recall that for every $\mu>0$, the function
$$w^\mu(x):=\frac1\mu w(\mu x)$$
is the solution of the BVP in the domain $\frac1\mu\Omega$. Let us list a few properties of the sequence $(w^\mu)_{\mu>0}$~:
\begin{itemize}
\item For $\epsilon<\eta$, one has the lower bound (maximum principle) $w^\epsilon\ge w^\eta$ in $\frac1\eta\Omega$.
\item For $\epsilon<\eta$, one has the Lipschitz estimate 
$$|\log w^\epsilon(x)-\log w^\epsilon(y)|\le d_H^\epsilon(x,y)\le d_H^\eta(x,y),\qquad\forall x,y\in\frac1\eta\Omega,$$
where we have denoted $d^\eta_H$ the Hilbert distance in $\frac1\eta\Omega$.
\item By the maximum principle, $w^\mu\le V$ for every $\mu>0$.
\end{itemize}
The Lipschitz estimates ensures that the PDE remains uniformly elliptic in every compact subdomain of the cone ${\cal C}$. Therefore the theory of elliptic regularity applies: every derivative $D^\beta w^\mu$ remains bounded as $\mu\rightarrow0^+$, on every compact subdomain of ${\cal C}$. The (monotonic) limit $w^0$ is again a solution of the PDE. In addition, it satisfies $(w^0)^\mu=w^0$, which means that it is homogeneous of degree one.
Because of the upper bound $w^0\le V$, we know that $w^0$ vanishes along the boundary. All this implies that $w^0$ is identical to $V$. 

Let us know select two points $x,y$ on the same ray $L$, close to the origin. The asymptotics above gives $|\log w(y)-\log w(x)|\sim|\log|y|-\log|x||\sim d_H(x,y)$. This implies $c_\Omega\ge1$. With Lemma \ref{l:LipH}, we conclude that $c_\Omega=1$.

\enpr

\bigskip

Another interesting situation is that of the equation 
$$\Delta w+\frac1w=0.$$
When $n=1$, and therefore the domain is $I=(0,\ell)$, the ODE can be integrated by hand and we find $u(x)\sim C_\ell x\sqrt{-2\log u}$. This quasi-Lipschitz behaviour at the boundary implies $c_I\ge1$, and therefore $c_I=1$.

\bigskip

The situation is significantly better for our fundamental example~:
\begin{prop}
Consider the BVP for the equation (\ref{eq:RP}). When $\Omega$ is a disk (hence $n=2$), we have $c_\Omega=\frac12\,.$
\end{prop} 

\bepr

By scaling, we may work in $D=D(0;1)$. Then
$$w(x)=\sqrt{\frac{1-|x|^2}2}\,,$$
a rare case where the solution is known in close form. In particular, the H\"older singularity of exponent $\frac12$ implies $c_\Omega\ge\frac12$\,.
On the other hand $d_H$ is given by
$$d_H(x,y)=2\log\left(1-x\cdot y+\sqrt{|y-x|^2-|x\wedge y|^2}\right)-\log(1-|x|^2)(1-|y|^2).$$
There remains the inequality $c_\Omega\le\frac12\,$.
The inequality $|\log w(x)-\log w(y)|\le \frac12d_H(x,y)$ to prove is equivalent to 
$$1-|x|^2\le1-x\cdot y+\sqrt{|y-x|^2-|x\wedge y|^2}.$$
It is implied by
$$(x\cdot(y-x))^2+|x\wedge y|^2\le|y-x|^2,$$
which is true because the left-hand side equals $(x\cdot(y-x))^2+|x\wedge (y-x)|^2=|x|^2|y-x|^2$, and on the other hand $|x|<1$. 

\enpr

\bigskip

We now show that the assumption made in Proposition \ref{p:coneL} is always met in our fundamental example. The cone ${\cal C}$ is a sector $S_\alpha$ of aperture $\alpha\in(0,\pi)$.
\begin{prop}
The BVP for the fundamental example (\ref{eq:RP}) is solvable in any planar sector $S_\alpha$.
\end{prop}
\begin{cor}
Let $\Omega$ be a planar open convex domain. Let us restrict to the equation (\ref{eq:RP}). If $\partial\Omega$ has a kink, then
$$c_\Omega=1.$$
\end{cor}

\bigskip

\bepr

Let us work in polar coordinates. The sector is
$$S_\alpha:=\{re^{i\theta}\,|\,\theta\in(0,\alpha)\}.$$
The self-similar solution is written $w_\alpha(x)=rA(\theta)$. The boundary condition is $A(0)=A(\alpha)=0$.

With $\nabla w_\alpha=A\vec e_r+A'\vec e_\theta$, the ODE satisfied by $\theta\mapsto A(\theta)$ is
$$\frac A{\sqrt{1+A^2+A'^2}}+\left(\frac{A'}{\sqrt{1+A^2+A'^2}}\right)'+\frac2{A\sqrt{1+A^2+A'^2}}=0,$$
that is
\begin{equation}\label{eq:edoasec}
A(1+A^2)(A''+A)+2(1+A^2+A'^2)=0.
\end{equation}
The solutions of (\ref{eq:edoasec}) may not be constant. 
This equation can therefore be integrated once, into
\begin{equation}\label{eq:edoapr}
A^4(1+A^2+A'^2)=C(1+A^2)^2,
\end{equation}
for some positive constant $C$. This autonomous ODE has the form $A'^2=F_C(A)$ where $F_C$ is positive over $(0,A^*)$ with
$$A^*=\sqrt{\frac12\left(C+\sqrt{C^2+4C}\right)\,}.$$ 
The Cauchy problem 
$$A'=\sqrt{F_C(A)}, \qquad A(0)=0$$
admits a unique maximal solution $A_C$ on an interval $[0,\ell]$, with
$$\ell=\int_0^{A^*}\frac{dA}{\sqrt{F_C(A)}}\,,$$
and we have $A_C'(\ell)=0$. The maximum principle tells us that at fixed $x$, the map $C\mapsto A_C(x)$ is increasing. In particular, $C\mapsto\ell$ is increasing; obviously, it is continuous too.

Let us compute the limits $\ell(0)$ and $\ell(+\infty)$. We have
$$\ell=\int_0^{A^*}\frac{A^2dA}{\sqrt{(1+A^2)(A^{*2}-A^2)(CA^{*-2}+A^2)}}=A^{*3}\int_0^1\frac{s^2ds}{\sqrt{(1+A^{*2}s^2)(1-s^2)(C+A^{*4}s^2)}}\,.$$
When $C\rightarrow0^+$, one has $A^{*2}\sim \sqrt C$ and therefore
$$\ell\sim C^{1/4}\int_0^1\frac{s^2ds}{\sqrt{1-s^4}}\rightarrow0,$$
whence $\ell(0)=0$. When instead $C\rightarrow+\infty$, we have $A^{*2}\sim C$ and
$$\ell=\int_0^1\frac{s^2ds}{\sqrt{(s^2+A^{*-2})(1-s^2)(s^2+CA^{*-4})}}\rightarrow\int_0^1\frac{ds}{\sqrt{1-s^2}}=\frac\pi2\,.$$

Extending $A_C$ by parity, we obtain a solution $A_C$ of (\ref{eq:edoapr}) vanishing at $0$ and $2\ell$, where $2\ell$ ranges from $0$ to $\pi$ when $C\in(0,+\infty)$. Therefore, there exists a unique $C$ for which $2\ell=\alpha$. Then $w_\alpha=rA_C(\theta)$ is the announced solution.

\enpr

\section{Existence proof}\label{s:exist}

So far, we have proved that if the solution $w$ of the BVP in $\Omega$ exists, then it enjoys a finite upper bound $w_+$, a positive lower bound $w_-$, and a Lipschitz estimate over every compact subdomain $K\subset\Omega$. This ensures that the linear operator $\sum_{i,j}a_{ij}(\nabla w)\partial_i\partial_j$ is uniformly elliptic on relatively compact subdomains. From regularity theory \cite{GT}, we deduce locally uniform estimates of  derivatives $\partial^\beta w$ of every order.

Our existence proof deals first with a modified problem, from which the singularity at the boundary has been removed, and the {\em a priori} uniform ellipticity has been restored.

\subsection{Relaxation of the boundary condition}

Let $\epsilon>0$ be given, we considered the BVP formed by the PDE (\ref{eq:PDE}), together with the boundary condition
\begin{equation}\label{eq:BCeps}
w=\epsilon\qquad\hbox{over }\partial\Omega.
\end{equation}
Because of the maximum principle, the solution $w_{\epsilon}$ must be unique and satisfy $w_\epsilon>\epsilon$ in $\Omega$. We may therefore replace the singularity $\frac1w$ in (\ref{eq:PDE}) by a smooth positive, decreasing function $g_\epsilon(w)$ which coincides with $\frac1w$ over $(\epsilon,+\infty)$. 

An upper barrier $w_{+,\epsilon}$ can be constructed by the following procedure. For every slab $\Pi_\nu$ containing $\Omega$, we consider the function
$$z_{\epsilon,\nu}(x)=\lambda W\left(\frac{x\cdot\nu-\alpha}\lambda\right)$$
with the same $W$ provided by Lemma \ref{l:W}. The parameters are chosen so that $z=\epsilon$ on the boundary of the slab:
$$\lambda=\frac{\alpha_+(\nu)-\alpha_-(\nu)}{1-2s},\qquad \alpha=\frac{1-s}{1-2s}\,\alpha_-(\nu)-\frac s{1-2s}\,\alpha_+(\nu),\qquad s=s(\epsilon):=W^{-1}(\epsilon).$$
Because $z_{\epsilon,\nu}$ is a solution in $\Pi_\nu$, it is a super-solution in $\Omega$. Therefore our upper barrier is 
$$w_{+,\epsilon}=\inf_{\nu\in S^{n-1}}z_{\epsilon,\nu}.$$
We point out that $|\nabla z_{\epsilon,\nu}|\le W'(s(\epsilon))$. This implies the same bound for $w_{+,\epsilon}$.

Let now $\bar a$ be a smooth numerical function that coincides with $a$ over $[0,W'(s(\epsilon))]$, such that $r\mapsto r\bar a(r)$ is increasing and $\bar a$ is constant over $[1+W'(s(\epsilon)),+\infty)$. Then the functions $z_{\epsilon,\nu}$ are super-solutions of the modified PDE
\begin{equation}
\label{eq:modPDE}
{\rm div}(\bar a(|\nabla w|)\nabla w)+\sigma g_\epsilon(w)F(|\nabla w|)=0
\end{equation}
as well, for every $\sigma\in[0,1]$ ; they are actually solutions when $\sigma=1$. The BVP (\ref{eq:modPDE},\ref{eq:BCeps}) admits therefore the upper barrier function $w_{+,\epsilon}$. Because a solution satisfies $\epsilon\le w\le w_{+,\epsilon}$, and since $w_{+,\epsilon}=\epsilon$ on the boundary, one infers that the normal derivative $\partial_\nu w$ at the boundary is bounded by that of $w_{+,\epsilon}$, that is by $W'(s(\epsilon))$. Then, because a PDE of the form above enjoys a maximum principle for derivatives, we find that for any solution of (\ref{eq:modPDE},\ref{eq:BCeps}), one has $|\nabla w|\le W'(s(\epsilon))$. 

All this, together with Theorem 11.3 of \cite{GT}, shows that the map $T$, defined by $w\mapsto z=Tw$ if
$${\rm div}(\bar a(|\nabla w|)\nabla z)+\sigma g_\epsilon(w)F(|\nabla w|)=0,\qquad z|_{\partial\Omega}=\epsilon$$
admits a fixed point $w_\epsilon$, which is a classical solution
of (\ref{eq:modPDE},\ref{eq:BCeps}). It satisfies the expected bounds
\begin{equation}\label{eq:espbd}
\epsilon\le w_\epsilon\le w_{+,\epsilon}.
\end{equation}
These bounds ensure that $w_\epsilon$ is actually a solution of (\ref{eq:PDE},\ref{eq:BCeps}). We point out that $w_\epsilon$ is unique.

\subsection{Passage to the limit}

We now prove that the $w_\epsilon$'s satisfy uniform estimates. On the one hand, the same rescaling as before can be used: if $p,q\in\Omega$, and if $\Omega-q\subset\lambda(\Omega-p)$, then
$$x\mapsto\lambda w_{\epsilon}\left(\frac x\lambda+p\right)$$
solves (\ref{eq:modPDE}) in $\Omega-q$, and is $\ge\lambda\epsilon\ge\epsilon$ over $\partial(\Omega-q)$. By the MP, we deduce
$$w_{\epsilon}(x+q)\le\lambda w_{\epsilon}\left(\frac x\lambda+p\right).$$
Setting $x=0$ in the inequality above, we obtain
$$w_{\epsilon}(q)\le\lambda w_{\epsilon}(p).$$
This shows that $w_{\epsilon}$ is Lipschitz with respect to the Hilbert metric, with Lipschitz constant $\le1$.

On the other hand, the same lower barrier $w_-$ applies to the modified BVP, and the upper barrier $w_{+,\epsilon}$ converges uniformly towards $w_+$ as $\epsilon\rightarrow0+$. By regularity theory, we therefore obtain uniform bounds for higher derivatives in every compact subdomain.

By Ascoli--Arzela and a diagonal procedure, we may extract from $(w_\epsilon)_{\epsilon>0}$ a subsequence that converges in ${\cal C}^{1,\beta}_{loc}(\Omega)$ for some $\beta>0$, to some limit function $w$. We may pass to the limit in (\ref{eq:PDE}), so that $w$ solves the PDE. On the other hand, passing to the limit in $w_-\le w_\epsilon\le w_{+,\epsilon}$ yields $w_-\le w\le w_+$. In particular, $w\in{\cal C}(\overline{\Omega})$ and $w$ satisfies the boundary condition (\ref{eq:DBC}). This ends the proof of Theorem \ref{th:main}.


\begin{thebibliography}{00}


 
   \bibitem{And} M. T. Anderson. Complete minimal varieties in hyperbolic space. {\em Inventiones mathematicae}, {\bf 69} (1982), pp 477--494.
 
    \bibitem{GT} D. Gilbarg, N. Trudinger. {\em Elliptic Partial Differential Equations of Second Order}. Classics in Mathematics. Springer-Verlag (2001), Heidelberg.

   \bibitem{Hil} D. Hilbert. Ueber die gerade Linie als k\"urzeste Verbindung zweier Punkte. {\em Mathematische Annalen}, {\bf 46} (1895), pp 91--96.
   
   \bibitem{Lin} Fang Hua Lin. On the Dirichlet problem for minimal graphs. {\em Inventiones mathematicae}, {\bf 96} (1989), pp 593--612.

   \bibitem{IdM} L. Marquis. G\'eom\'etrie de Hilbert. {\em Images des Math\'ematiques}, CNRS (2015). \\ {\tt  http://images.math.cnrs.fr/Geometrie-de-Hilbert.html}.
 
   \bibitem{Ser_Chap} D. Serre. Multi-dimensional shock interaction for a Chaplygin gas. {\em Arch. Rational Mech. Anal.}, {\bf191} (2009), pp 539--577.
        
   \bibitem{Ser_log} D. Serre. Gradient estimate in terms of a Hilbert-like distance, for minimal surfaces and Chaplygin gas. {\em Comm. Partial Diff. Equ.},  {\bf41} (2016), pp 774--784.
        
   \bibitem{Wal} C. Walsh. Gauge-reversing maps on cones , and Hilbert and Thompson isometries.  Preprint arXiv:1312.7871 [math.MG] (december 2013).


\end{thebibliography}
\end{document}